\documentclass{amsart}
\usepackage{wasowicz}

\newcommand{\K}{\mathcal{K}}
\renewcommand{\le}{\leqslant}
\renewcommand{\ge}{\geqslant}
\newcommand{\norm}[1]{\|#1\|_{\infty}}

\info{Published: Math. Inequal. Appl. 13~(2010), 165--174.}

\begin{document}

\title{On some extremalities in the approximate integration}
\author{\SW}
\address{\SWaddr}
\email{\SWmail}
\date{February 2008}
\begin{abstract}
 Some extremalities for quadrature operators are proved for convex functions of
 higher order. Such results are known in the numerical analysis, however they are
 often proved under suitable differentiability assumptions. In our considerations
 we do not use any other assumptions apart from higher order convexity itself.
 The obtained inequalities refine the inequalities of Hadamard type. They are
 applied to give error bounds of quadrature operators under the assumptions weaker
 from the commonly used.
\end{abstract}
\maketitle

\section{Introduction}

In the theory of convex functions the Hermite--Hadamard inequality
\begin{equation}\label{HH}
 f\left(\frac{a+b}{2}\right)\le\frac{1}{b-a}\int_a^bf(x)dx\le\frac{f(a)+f(b)}{2},
\end{equation}
which holds for convex functions (and, in fact, characterizes them), plays a
very important role. The first inequality follows by the existence of a support
line for $f$ at the midpoint, while the second one can be obtained using the
fundamental property of convexity stating that a graph of a convex function $f$
lies on $[a,b]$ below the chord joining the points $\bigl(a,f(a)\bigr)$,
$\bigl(b,f(b)\bigr)$. We have also the following
\begin{obs}\label{mid_trap}
 If a real function $f$ is convex on an interval $[a,b]$ then
 \begin{equation}\label{ineq:mid_trap}
  f\left(\frac{a+b}{2}\right)\le\sum_{i=1}^N\lambda_if(\xi_i)\le\frac{f(a)+f(b)}{2}
 \end{equation}
 for any $N\in\N$, $\xi_1,\dots,\xi_N\in[a,b]$ and
 $\lambda_1,\dots,\lambda_N\ge 0$ with $\sum_{i=1}^N\lambda_i=1$ such that
 $\sum_{i=1}^N\lambda_i\xi_i=\frac{a+b}{2}$.
\end{obs}
\begin{proof}
 The first inequality is an immediate consequence of convexity, the second
 one we prove similarly to the second inequality of~\eqref{HH}.
\end{proof}

The term on the left hand side of~\eqref{ineq:mid_trap} is connected with the
midpoint rule of the approximate integration, while the term on the right hand
side is connected with the trapezoidal rule. Then the
inequality~\eqref{ineq:mid_trap} can be regarded as an example of an
extremality for quadrature operators. Many extremalities are known in the
numerical analysis (cf. \cite{BraSch81}, cf. also \cite{BraPet}
and the references therein). The numerical analysts prove them using the
suitable differentiability assumptions. As we will show in this paper, for
convex functions of higher order some extremalities can be obtained without
assumptions of this kind, using only the higher order convexity itself. The
suppor--type properties play here the crucial role. A~general theorem of this
nature was recently proved by the author in~\cite{Was07JMAA}. The obtained
extremalities are useful in proving error bounds of quadrature operators under
regularity assumptions weaker from the commonly used. The results of this sort
are also known, however our method seems to be quite easy. But the price we
must pay is high: the obtained error bounds are far to be optimal (cf.
\cite{BraPet}). Some results concerning the inequalities between the
quadrature operators and error bounds of quadrature rules, which are partial
cases of our results, can be found in author's earlier papers
\cite{Was07JIPAM,Was06JIPAM,Was08JIPAM2}. The paper \cite{Was08MIA}
contains the extension of our results concerning convex functions to the
functions of several variables.

For $n\in\N$ denote by $\Pi_n$ the space of all polynomials of degree at
most~$n$. Recall that a~linear functional $\T$ defined on a linear space $X$ of
(not necessarily all) functions mapping some nonempty set into $\R$ is called
\emph{positive} if
\[
 f\le g\implies\T(f)\le\T(g)
\]
for any $f,g\in X$. An important class of positive linear operators form the
conical combinations of the involved function  at appropriately chosen points of a
domain. Obviously, if a domain is a real interval, then quadrature operators
with nonnegative coefficients are linear and positive.

Dealing with a problem of approximate computation of the integral over an
interval $[a,b]$ it is enough to change the variable and to compute it over a
fixed interval. The interval $[-1,1]$ is frequently used. For a Riemann
integrable function $f:[-1,1]\to\R$ let
\[
 \I(f):=\int_{-1}^1 f(x)dx.
\]
\begin{defn}
 Let $\T$ be a linear functional defined on a linear space of (not necessarily
 all) Riemann integrable functions mapping $[-1,1]$ into $\R$ containing $\Pi_n$.
 We say that $\T$ is \emph{exact on} $\Pi_n$ if $\T(p)=\I(p)$ for all
 $p\in\Pi_n$.
\end{defn}

\section{Convex functions of higher order}

Recall that the divided differences are defined as follows: $[x_1;f]:=f(x_1)$
and for $k\in\N$
\[
 [x_1,\dots,x_{k+1};f]:=\frac{[x_2,\dots,x_{k+1};f]-[x_1,\dots,x_k;f]}{x_{k+1}-x_1}.
\]

\begin{defn}\label{nconv_def}
 Let $I\subset\R$ be an interval and $n\in\N$. A function $f:I\to\R$ is
 $n$-\emph{convex} if $[x_1,\dots,x_{n+2};f]\ge 0$ for any distinct
 $x_1,\dots,x_{n+2}\in I$.
\end{defn}

There is an easy to imagine geometrical equivalent condition of $n$-convexity
(for the proof cf. e.g. \cite{Kuc85,Pop34}).
\begin{prop}
 A function $f$ is $n$-convex if and only if for any $n+1$ distinct points
 $x_1,\dots,x_{n+1}\in I$, the graph of the (unique) polynomial $p\in\Pi_n$
 interpolating $f$ at these points, passing through each point
 $\bigl(x_i,f(x_i)\bigr)$, $i=1,\dots,n+1$, changes the side of the graph of $f$
 (always $p\ge f$ on $[x_n,x_{n+1}]$).
\end{prop}
Then trivially 1-convexity reduces to the classical convexity.

Convex functions of higher order are very well known and investigated (see e.g.
\cite{Bul71,Kuc85,PinWul05,Pop34,RobVar73}). But up to now
there is no common terminology, which sometimes may be confusing.

The first person who dealt with the topic in question was Hopf. He considered
in his dissertation~\cite{Hop26} from 1926 the functions with nonnegative
divided differences without naming them at all. The notion of higher order
convexity was introduced by Popoviciu in his famous dissertation
\cite{Pop34} from 1934 exactly in the sense of the above
Definition~\ref{nconv_def}. The Kuczma's monograph \cite{Kuc85} devoted to
functional equations and inequalities in several variables as well as the
classical Roberts and Varberg's book on convex functions \cite{RobVar73}
use the same terminology (according to which an ordinary convex function is
1-convex). During the years another way of naming convex functions of higher
order became popular. Some authors (cf. e.g. \cite{Bul71,PinWul05}) call
a function $f$ to be $n$-convex if $[x_1,\dots,x_{n+1};f]\ge 0$ (then a convex
function is 2-convex). Now these two terminologies appear simultaneously in
the literature. The first one is concentrated on the maximal degree of the
interpolating polynomial, while the accent of the second one is put on the
dimension of the space of polynomials of degree not exceeding some natural
number. Then the second terminology is more coherent with convexity with
respect to Chebyshev systems (cf. \cite{KarStu66}). Both conventions have
some advantages and disadvantages and it is not the author's intention to judge
neither which one is better nor which one is classical.

Having in mind the above remarks let us declare that in this paper we
understand the higher order convexity in the sense of
Definition~\ref{nconv_def}.

Convex functions of higher order have many regularity properties. For details
see \cite{Kuc85,Pop34,RobVar73}. The paper \cite{PinWul05}
contains a brief survey of the topic given in one place. Below we list the
properties which either we use in the paper or we discuss below.

\begin{thm}
 If $f:[a,b]\to\R$ is $n$-convex then $f$ is continuous on $(a,b)$ and bounded
 on $[a,b]$.
\end{thm}

\begin{cor}
 If $f:[a,b]\to\R$ is $n$-convex then $f$ is Riemann integrable.
\end{cor}

\begin{thm}\label{nconv_diff}
 The real function $f$ defined on an open interval $I$ is $n$-convex if and only
 if $f^{(n-1)}$ is convex on $I$.
\end{thm}

\begin{cor}
 If the real function $f$ defined on an open interval $I$ is $n$-convex
 then $f^{(n)}_-$, $f^{(n)}_+$ exist on I and $f^{(n)}$ exists almost
 everywhere on $I$.
\end{cor}

Notice that there are $n$-convex functions which are not $n$ times
differentiable (e.g. $f(x)=|x|$ for $n=1$). Sometimes what is proved under
differentiability assumptions, holds in fact for any $n$-convex function
without further assumptions. We return to this matter in
Section~\ref{extremalities}. However, the following result requiring the
differentiability assumption seems to be important (cf.
\cite{Kuc85,Pop34,RobVar73}, for a quick reference cf. also
\cite[Theorems A and B]{Was06JIPAM} and \cite[Theorem D]{Was07JMAA}).

\begin{thm}\label{nconv_diff_nonnegative}
 Assume that $f:[a,b]\to\R$ is ($n+1$)-times differentiable on $(a,b)$ and
 continuous on $[a,b]$. Then~$f$ is $n$-convex if and only if
 $f^{(n+1)}(x)\ge 0$, $x\in(a,b)$.
\end{thm}

It is well known that a convex function defined on a real interval admits an
affine support at every interior point of a domain. In the paper
\cite{Was07JMAA} we have proved a general suppor--type result for convex
functions of higher order. Four special cases
(\cite[Corollaries~8-11]{Was07JMAA}) play the crucial role in the proofs
presented in this paper.

\begin{thm}\label{supp_G}
 If $f:[a,b]\to\R$ is $(2n-1)$-convex and
 $x_1,\dots,x_n\in (a,b)$, then there exists a~polynomial
 $p\in\Pi_{2n-1}$ such that $p(x_i)=f(x_i)$, $i=1,\dots,n$,
 and $p\le f$ on $[a,b]$.
\end{thm}
\begin{thm}\label{supp_L}
 If $f:[a,b]\to\R$ is $(2n-1)$-convex and $x_1=a$,
 $x_2,\dots,x_n\in (a,b)$, $x_{n+1}=b$, then there exists a~polynomial
 $p\in\Pi_{2n-1}$ such that $p(x_i)=f(x_i)$, $i=1,\dots,n+1$, and
 $p\ge f$ on $[a,b]$.
\end{thm}
\begin{thm}\label{supp_Rl}
 If $f:[a,b]\to\R$ is $2n$-convex, $x_1=a$,
 $x_2,\dots,x_{n+1}\in (a,b)$, then there exists a~polynomial
 $p\in\Pi_{2n}$ such that $p(x_i)=f(x_i)$, $i=1,\dots,n+1$,
 and $p\le f$ on $[a,b]$.
\end{thm}
\begin{thm}\label{supp_Rr}
 If $f:[a,b]\to\R$ is $2n$-convex, $x_1,\dots,x_n\in
 (a,b)$ and $x_{n+1}=b$, then there exists a~polynomial $p\in\Pi_{2n}$
 such that $p(x_i)=f(x_i)$, $i=1,\dots,n+1$, and $p\ge f$ on $[a,b]$.
\end{thm}

\section{Extremalities for quadrature operators}\label{extremalities}

What we recall below is very well known from the numerical analysis (cf. e.g.
\cite{BesPal02,Hil56,Ral65,Szeg75,WeiGL,WeiL,WeiR}). Let $P_n$ be the $n$-th
degree member of the sequence of Legendre polynomials.

\subsection*{Gaus--Legendre quadratures}

For $f:[-1,1]\to\R$ and $n\in\N$ let
\[
 \G_n(f):=\sum_{i=1}^nw_if(x_i),
\]
where $x_1,\dots,x_n$ are the roots of $P_n$ (which are real, distinct and
belong to $(-1,1)$) and
\[
 w_i=\frac{2(1-x_i^2)}{(n+1)^2P_{n+1}^2(x_i)},\quad i=1,\dots,n.
\]
Then $\G_n$ is exact on $\Pi_{2n-1}$. If $f\in\C^{2n}([-1,1])$ then
\begin{equation}\label{err:G}
 \I(f)=\G_n(f)+\frac{2^{2n+1}(n!)^4}{(2n+1)[(2n)!]^3}f^{(2n)}(\xi)
\end{equation}
for some $\xi\in(-1,1)$.

\subsection*{Lobatto quadratures}

For $f:[-1,1]\to\R$ let $\Lob_2(f):=f(-1)+f(1)$. For $n\in\N$, $n\ge 3$, let
\[
 \Lob_n(f):=w_1f(-1)+w_nf(1)+\sum_{i=2}^{n-1}w_if(x_i),
\]
where $x_2,\dots,x_{n-1}$ are the roots of $P_{n-1}'$ (which are also real,
distinct and belong to $(-1,1)$) and
\[
 w_1=w_n=\frac{2}{n(n-1)},\quad w_i=\frac{2}{n(n-1)P_{n-1}^2(x_i)},
 \quad i=2,\dots,n-1.
\]
Then $\Lob_n$ is exact on $\Pi_{2n-3}$. If $f\in\C^{2n-2}([-1,1])$ then
\begin{equation}\label{err:L}
 \I(f)=\Lob_n(f)
      -\frac{n(n-1)^32^{2n-1}[(n-2)!]^4}{(2n-1)[(2n-2)!]^3}f^{(2n-2)}(\xi)
\end{equation}
for some $\xi\in(-1,1)$.

\subsection*{Radau quadratures}

For $f:[-1,1]\to\R$ and $n\in\N$, $n\ge 2$, let
\[
 \Rad_n^l(f):=w_1f(-1)+\sum_{i=2}^nw_if(x_i),
\]
where $x_2,\dots,x_n$ are the roots of the polynomial
\[
 Q_{n-1}(x)=\frac{P_{n-1}(x)+P_n(x)}{x+1}
\]
(again real, distinct and belonging to $(-1,1)$) and
\[
 w_1=\frac{2}{n^2},\quad w_i=\frac{1}{(1-x_i)[P_{n-1}'(x_i)]^2},
 \quad i=2,\dots,n.
\]
Then $\Rad_n^l$ is exact on $\Pi_{2n-2}$. If $f\in\C^{2n-1}([-1,1])$ then
\begin{equation}\label{err:R}
 \I(f)=\Rad_n^l(f)
      +\frac{2^{2n-1}n[(n-1)!]^4}{[(2n-1)!]^3}f^{(2n-1)}(\xi)
\end{equation}
for some $\xi\in(-1,1)$.

In \cite{Was07JMAA} we also considered the operator
\[
 \Rad_n^r(f):=\Rad_n^l\bigl(f(-\,\cdot)\bigr).
\]
It was, in fact, defined in terms of
orthogonal polynomials. However, these two definitions coincide. This is not
difficult to check. We would not like to go into details since this is not the
goal of the paper. Let us only mention that (changing the way of naming and
numbering the abscissas and weights) we have
\[
 \Rad_n^r(f)=\sum_{i=1}^nw_if(y_i)+w_nf(1)
\]
and $\Rad_n^r$ is exact on $\Pi_{2n-2}$. The error term of $\Rad_n^r$ is
similar to \eqref{err:R}, precisely, if $f\in\C^{2n-1}([-1,1])$ then
\begin{equation}\label{err:Rr}
 \I(f)=\Rad_n^r(f)
      -\frac{2^{2n-1}n[(n-1)!]^4}{[(2n-1)!]^3}f^{(2n-1)}(\eta)
\end{equation}
for some $\eta\in(-1,1)$.

As we can see, all the weights of the above quadratures are positive, so these
operators are positive. This is also the case for many other quadratures.
However, there are the quadratures with negative coefficients (e.g. among the
Newton--Cotes formulas).

Now we can prove the main results of this section.

\begin{thm}\label{Gauss_Lobatto}
 Fix $n\in\N$. Let $\T$ be the positive linear operator defined (at least) on a~linear
 subspace of $\R^{[-1,1]}$ generated by a cone of $(2n-1)$-convex functions.
 Assume that $\T$ is exact on $\Pi_{2n-1}$. If a function $f:[-1,1]\to\R$ is
 $(2n-1)$-convex then
 \begin{equation}\label{ineq:Gauss_Lobatto}
  \G_n(f)\le\T(f)\le\Lob_{n+1}(f).
 \end{equation}
\end{thm}
\begin{proof}
 By Theorems~\ref{supp_G}~and~\ref{supp_L} there exist two polynomials
 $p,q\in\Pi_{2n-1}$ interpolating $f$ at the abscissas of the operators $\G_n$,
 $\Lob_{n+1}$, respectively, such that $p\le f\le q$ on $[-1,1]$. Since
 $\G_n=\Lob_{n+1}=\I$ on $\Pi_{2n-1}$, we get
 \[
  \G_n(f)=\G_n(p)=\I(p)=\T(p)\le\T(f)\le\T(q)=\I(q)=\Lob_{n+1}(q)=\Lob_{n+1}(f).
  \qedhere
 \]
\end{proof}

\begin{thm}\label{Radau}
 Fix $n\in\N$. Let $\T$ be the positive linear operator defined (at least) on a~linear
 subspace of $\R^{[-1,1]}$ generated by a cone of $2n$-convex functions.
 Assume that $\T$ is exact on $\Pi_{2n}$. If a function $f:[-1,1]\to\R$ is
 $2n$-convex then
 \begin{equation}\label{ineq:Radau}
  \Rad_{n+1}^l(f)\le\T(f)\le\Rad_{n+1}^r(f).
 \end{equation}
\end{thm}
\begin{proof}
 Use Theorems~\ref{supp_Rl}~and~\ref{supp_Rr} and the abscissas of the operators
 $\Rad_{n+1}^l$, $\Rad_{n+1}^r$, respectively, and argue similarly as in the
 proof of Theorem~\ref{Gauss_Lobatto}.
\end{proof}

We would like to emphasize two particular cases of the above results. The first
one concerns the inequalities of Hadamard type. The assertions of the Corollary
below were proved in \cite[Propositions 12 and 13]{Was07JMAA} (cf. also the
earlier paper \cite{BesPal02}, where these results were obtained by another method).

\begin{cor}
 Fix $n\in\N$. If $f:[-1,1]\to\R$ is $(2n-1)$-convex then
 $\G_n(f)\le\I(f)\le\Lob_{n+1}(f)$. If $f:[-1,1]\to\R$ is $2n$-convex then
 $\Rad_{n+1}^l(f)\le\I(f)\le\Rad_{n+1}^r(f)$.
\end{cor}

\begin{proof}
 Use Theorems \ref{Gauss_Lobatto} and \ref{Radau} for $\T=\I$.
\end{proof}

The second important case is connected with quadrature operators.

\begin{cor}\label{quadrature}
 Fix $n,N\in\N$, $\xi_1,\dots,\xi_N\in[-1,1]$ and $\lambda_1,\dots,\lambda_N\ge 0$.
 Let
 \[
  \T(f):=\sum_{i=1}^N\lambda_if(\xi_i)\quad\text{for }f:[-1,1]\to\R.
 \]
 \begin{enumerate}[\upshape(i)]
  \item\label{quadrature:i}
   If $\T$ is exact on $\Pi_{2n-1}$, then $\G_n(f)\le\T(f)\le\Lob_{n+1}(f)$
   for any $(2n-1)$-convex function $f:[-1,1]\to\R$.
  \item\label{quadrature:ii}
   If $\T$ is exact on $\Pi_{2n}$, then $\Rad_{n+1}^l(f)\le\T(f)\le\Rad_{n+1}^r(f)$
   for any $2n$-convex function $f:[-1,1]\to\R$.
 \end{enumerate}
\end{cor}

\begin{proof}
 The operator $\T$ trivially fulfils the assumptions of Theorems
 \ref{Gauss_Lobatto} and~\ref{Radau}, respectively.
\end{proof}

In the numerical analysis the inequalities the above type are called
\emph{extremalities}. The extremalities of Corollary~\ref{quadrature}~\itemref{quadrature:i}
were earlier proved in \cite[Theorem 6]{BraSch81} under the assumption
of $2n$-times differentiability. The proof given there is based on taking
double nodes. The author independently used in \cite{Was07JMAA} exactly the same
idea to prove support--type results of Corollaries 8-11 (quoted here in
Theorems \ref{supp_G}-\ref{supp_Rr}) with no use of any differentiability
assumptions. Thus, as we can see from the proof of Theorems~\ref{Gauss_Lobatto}
and \ref{Radau}, the extremalities in question are proved with no further
assumptions, except higher order convexity itself. So, our results are more
general than these of \cite{BraSch81}.

We underline that the inequalities of Corollary~\ref{quadrature} do not hold
for any quadrature operator $\T$. The exactness assumption (i.e. $\T=\I$ for
polynomials of the appropriate degree) is essential.

\begin{exmp}
  We have
  \begin{align*}
    \G_2(f)&=f\biggl(-\frac{\sqrt{3}}{3}\biggr)
            +f\biggl(\frac{\sqrt{3}}{3}\biggr),\\
    \G_3(f)&=\frac{8}{9}f(0)
            +\frac{5}{9}\biggl[
                          f\biggl(-\frac{\sqrt{15}}{5}\biggr)
                         +f\biggl(\frac{\sqrt{15}}{5}\biggr)
                        \biggr],\\
   \Rad_3^l(f)&=\frac{2}{9}f(-1)
               +\frac{16+\sqrt{6}}{18}f\biggl(\frac{1-\sqrt{6}}{5}\biggr)
               +\frac{16-\sqrt{6}}{18}f\biggl(\frac{1+\sqrt{6}}{5}\biggr)
  \end{align*}
  (cf. e.g. \cite{WeiGL,WeiR}), whence $\G_3(\exp)>\G_2(\exp)$ and
  $\Rad_3^l(\exp)>\G_2(\exp)$. The exponential function is convex of any order
  (cf. Theorem~\ref{nconv_diff_nonnegative}). Let $\T=\G_2$. Then the inequality of
  Corollary~\ref{quadrature}~\itemref{quadrature:i} does not hold for $n=3$, and that of \itemref{quadrature:ii}
  is not true for $n=2$. Notice that for $p(x)=x^4$ we have $\G_2(p)\ne\I(p)$,
  so $\G_2\ne\I$ both on $\Pi_5$ and on $\Pi_4$.
\end{exmp}

Now the question arises if there are other, i.e. non--quadrature, operators
approximating the integral, for which Theorems~\ref{Gauss_Lobatto} and
\ref{Radau} are applicable. The positive answer given below shows that the
extremalities for quadrature operators are special cases of more general
inequalities. Not the form of the operators considered (linear combination,
integral and so on) is important but two things play the key role: positiveness
and exactness for polynomials of the appropriate degree.

\begin{exmp}
  For a Riemann--integrable function $f:[-1,1]\to\R$ let
  \[
   \T(f):=\frac{3}{11}[f(-1)+f(1)]
         +\frac{16}{11}\int_{-\frac{1}{2}}^{\frac{1}{2}}f(t)dt.
  \]
  Then $\T$ is a positive linear operator exact on $\Pi_3$. Using
  Theorem~\ref{Gauss_Lobatto} for $n=2$ we obtain
  $
   \G_2(f)\le\T(f)\le\Lob_3(f)
  $
  for a 3-convex function $f:[-1,1]\to\R$.
\end{exmp}

For the other non--quadrature operators approximating the integral cf. e.g.
\cite{BojPet01}.

\section{Error bounds of quadrature operators}

In this section we show that the extremalities of Corollary~\ref{quadrature}
may be applied to obtain the error bounds of the involved quadrature operator
$\T$ using the regularity assumptions weaker from the commonly used. The
results of this type are known in the numerical analysis. We would like to
point that the inequalities of Hadamard type may be used in the approximate
integration. But the results obtained by our method deliver error bounds which
are far to be optimal. This is a price we have to pay for simplicity. Error
bounds obtained in the numerical analysis under assumptions used by us are much
better (see \cite{BraPet} and the references therein).

For $f\in\C([-1,1])$ denote
$\norm{f}:=\sup\left\{\bigl|f(x)\bigr|\,:\,x\in[-1,1]\right\}$.

\begin{lem}\label{lem_eb}
 Fix $k\in\N$, $k\ge 2$. Let $\K$, $\T$ be linear operators defined on
 a linear subspace of\/ $\R^{[-1,1]}$ containing all the functions involved below
 with the following properties:
 \begin{enumerate}[\upshape(i)]
  \item\label{lem_eb:i}
   there exists an $\alpha>0$ such that $\I(f)\le\K(f)+\alpha\norm{f^{(k)}}$ for
   all $f\in\C^k([-1,1])$;
  \item\label{lem_eb:ii}
   If $f:[-1,1]\to\R$ is ($k-1$)-convex then $\K(f)\le\T(f)$;
  \item\label{lem_eb:iii}
   $\T(p)=\I(p)$ for $p(x)=x^k$.
 \end{enumerate}
 Then $\bigl|\I(f)-\T(f)\bigr|\le 2\alpha\norm{f^{(k)}}$ for any
 $f\in\C^k([-1,1])$.
\end{lem}
\begin{proof}
 By \eqref{lem_eb:i} and \eqref{lem_eb:ii} we get
 \begin{equation}\label{lem_eb_pf1}
  \I(f)-\T(f)\le\alpha\norm{f^{(k)}}
 \end{equation}
 for any ($k-1$)-convex function $f\in\C^k([-1,1])$.

 For an arbitrary function $f\in\C^k([-1,1])$ define now
 $g(x):=\frac{\norm{f^{(k)}}}{k!}x^k$. Then $g^{(k)}=\norm{f^{(k)}}$, whence
 $|f^{(k)}|\le g^{(k)}$ on $[-1,1]$, which implies
 $(g-f)^{(k)}\ge 0$ and $(g+f)^{(k)}\ge 0$ on $[-1,1]$. Therefore
 by Theorem~\ref{nconv_diff_nonnegative} the functions $g-f$, $g+f$ are
 ($k-1$)-convex on $[-1,1]$. By the triangle inequality
 \[
  \norm{(g-f)^{(k)}}\le 2\norm{f^{(k)}}\quad\text{and}\quad
  \norm{(g+f)^{(k)}}\le 2\norm{f^{(k)}}.
 \]
 Now we apply \eqref{lem_eb_pf1} to $g-f$ and $g+f$. Then the desired inequality
 follows by linearity, the assumption~\eqref{lem_eb:iii}  and the above inequalities.
\end{proof}

For $n\in\N$ let
\[
 \alpha_{2n}:=\frac{4^{n+1}(n!)^4}{(2n+1)[(2n)!]^3},\qquad
 \alpha_{2n+1}:=\frac{4^{n+1}(n+1)(n!)^4}{[(2n+1)!]^3}
\]

\begin{thm}\label{eb}
 Fix $k\in\N$, $k\ge 2$. Let $\T$ be a positive linear operator defined on a domain
 as in Lemma~\ref{lem_eb}. If $\T$ is exact on $\Pi_k$, then
 $\bigl|\I(f)-\T(f)\bigr|\le\alpha_k\norm{f^{(k)}}$ for any $f\in\C^k([-1,1])$.
\end{thm}

\begin{proof}
 If $k$ is even, $k=2n$, then use Lemma~\ref{lem_eb} for
 $\alpha=\frac{\alpha_{2n}}{2}$ and $\K=\G_n$. The condition $(i)$ is fulfilled
 by \eqref{err:G}, $(ii)$ holds by Theorem~\ref{Gauss_Lobatto} and $(iii)$ by
 the assumption.

 Similarly, if $k$ is odd, $k=2n+1$, then use Lemma~\ref{lem_eb} for
 $\alpha=\frac{\alpha_{2n+1}}{2}$ and $\K=\Rad_{n+1}^l$. The condition $(i)$ is
 fulfilled by \eqref{err:R}, $(ii)$ holds by Theorem~\ref{Radau} and
 $(iii)$ by the assumption.
\end{proof}

In the assertion of Theorem~\ref{Gauss_Lobatto} there is an operator
$\Lob_{n+1}$ on the right hand side of the inequality
\eqref{ineq:Gauss_Lobatto} and in the statement of Theorem~\ref{Radau} there
is an operator $\Rad_{n+1}^r$ on the right hand side of the inequality
\eqref{ineq:Radau}. We could prove the result similar to Theorem~\ref{eb}
involving (in the proof) these operators. However, the error bound obtained
in this way will not improve that of Theorem~\ref{eb}. Namely, the absolute
value of the constant of \eqref{err:L} (take $n+1$ instead of $n$) is greater
from the similar constant of \eqref{err:G}. For the operators $\Rad_{n+1}^l$
and $\Rad_{n+1}^r$ the absolute values of both constants of \eqref{err:R},
\eqref{err:Rr} are the same.

For the quadrature operators we immediately derive from Theorem~\ref{eb} the
following

\begin{cor}\label{eb:quadratures}
 Fix $k,N\in\N$, $k\ge 2$, $\xi_1,\dots,\xi_N\in[-1,1]$ and
 $\lambda_1,\dots,\lambda_N\ge 0$. Let
 \[
  \T(f):=\sum_{i=1}^N\lambda_if(\xi_i)\quad\text{for }f:[-1,1]\to\R.
 \]
 If $\T$ is exact on $\Pi_k$, then
 $\bigl|\I(f)-\T(f)\bigr|\le\alpha_k\norm{f^{(k)}}$
 for any $f\in\C^k([-1,1])$.
\end{cor}

Using this result we will now give the error bounds of Gauss--Legendre, Lobatto
and Radau quadratures under regularity assumptions weaker from the commonly
used. Denote by $\lfloor\cdot\rfloor$ the floor function, i.e.
$\lfloor x\rfloor=\max\{k\in\Z:k\le x\}$, $x\in\R$.

\begin{prop}
  Let $k,N\in\N$, $k\ge 2$ and $N>\left\lfloor\frac{k}{2}\right\rfloor$.
  If $f\in\C^k([-1,1])$ then
  $\bigl|\I(f)-\G_N(f)\bigr|\le\alpha_k\norm{f^{(k)}}$.
\end{prop}

\begin{proof}
  If $N>\left\lfloor\frac{k}{2}\right\rfloor$ then $2N-1\ge k$, whence
  $\G_N$ is exact on $\Pi_k$ (cf. \eqref{err:G}). Then the result follows
  immediately by Corollary~\ref{eb:quadratures}.
\end{proof}

\begin{prop}
  Let $k,N\in\N$, $k\ge 2$ and $N>\left\lfloor\frac{k}{2}\right\rfloor+1$.
  If $f\in\C^k([-1,1])$ then
  $\bigl|\I(f)-\Lob_N(f)\bigr|\le\alpha_k\norm{f^{(k)}}$.
\end{prop}

\begin{proof}
  This is also an immediate consequence of Corollary~\ref{eb:quadratures} since
  $N>\left\lfloor\frac{k}{2}\right\rfloor+1$ implies $2N-3\ge k$ and then
  $\Lob_N$ is exact on $\Pi_k$ (cf. \eqref{err:L}).
\end{proof}

\begin{prop}
  Let $k,N\in\N$, $k\ge 2$ and $N>\left\lfloor\frac{k+1}{2}\right\rfloor$. If
  $f\in\C^k([-1,1])$ then $\bigl|\I(f)-\Rad_N^l(f)\bigr|\le\alpha_k\norm{f^{(k)}}$.
  The same assertion holds for the operator $\Rad_N^r$.
\end{prop}

\begin{proof}
  By $N>\left\lfloor\frac{k+1}{2}\right\rfloor$ we get $2N-2\ge k$ and we can
  see (cf. \eqref{err:R}, \eqref{err:Rr}) that $\Rad_N^l$ and $\Rad_N^r$ are exact
  on $\Pi_k$, which, together with Corollary~\ref{eb:quadratures}, concludes the
  proof.
\end{proof}

\bibliographystyle{amsplain}
\bibliography{was_pub,was_own}

\providecommand{\bysame}{\leavevmode\hbox to3em{\hrulefill}\thinspace}
\providecommand{\MR}{\relax\ifhmode\unskip\space\fi MR }
\providecommand{\MRhref}[2]{%
  \href{http://www.ams.org/mathscinet-getitem?mr=#1}{#2}
}
\providecommand{\href}[2]{#2}
\begin{thebibliography}{10}

\bibitem{BesPal02}
Mihály Bessenyei and Zsolt Páles, \emph{Higher-order generalizations of
  {H}adamard's inequality}, Publ. Math. Debrecen \textbf{61} (2002), no.~3-4,
  623--643. \MR{1943721 (2003k:26021)}

\bibitem{BojPet01}
Borislav~D. Bojanov and Petar~P. Petrov, \emph{Gaussian interval quadrature
  formula}, Numer. Math. \textbf{87} (2001), no.~4, 625--643. \MR{1815728
  (2001k:65043)}

\bibitem{BraPet}
Helmut Brass and Knut Petras, \emph{Quadrature theory}, \\online:
  \texttt{http://www-public.tu-bs.de:8080/$\sim$petras/book/qt2.ps}.

\bibitem{BraSch81}
Helmut Brass and Gerhard Schmeisser, \emph{Error estimates for interpolatory
  quadrature formulae}, Numer. Math. \textbf{37} (1981), no.~3, 371--386.
  \MR{627111 (82j:65012)}

\bibitem{Bul71}
Peter Bullen, \emph{A criterion for {$n$}-convexity}, Pacific J. Math.
  \textbf{36} (1971), 81--98. \MR{0274681 (43 \#443)}

\bibitem{Hil56}
F.~B. Hildebrand, \emph{Introduction to numerical analysis}, McGraw-Hill Book
  Company, Inc., New York-Toronto-London, 1956. \MR{0075670 (17,788d)}

\bibitem{Hop26}
Eberhard Hopf, \emph{Über die {Z}usammenhänge zwischen gewissen höheren
  {D}ifferenzenquotienten reeller {F}unktionen einer reellen {V}ariablen und
  deren {D}ifferenzierbarkeitseigenschaften}, Ph.{D}. dissertation,
  Friedrich-Wilhelms-Universität, Berlin, 1926.

\bibitem{KarStu66}
Samuel Karlin and William~J. Studden, \emph{Tchebycheff systems: {W}ith
  applications in analysis and statistics}, Pure and Applied Mathematics, Vol.
  XV, Interscience publishers John Wiley \& Sons, New York-London-Sydney, 1966.
  \MR{0204922 (34 \#4757)}

\bibitem{Kuc85}
Marek Kuczma, \emph{An {I}ntroduction to the {T}heory of {F}unctional
  {E}quations and {I}nequalities}, Prace Naukowe Uniwersytetu Śląskiego w
  Katowicach, vol. 489, Państwowe Wydawnictwo Naukowe---Uniwersytet Śląski,
  Warszawa--Kraków--Katowice, 1985. \MR{788497 (86i:39008)}

\bibitem{PinWul05}
Allan Pinkus and Dan Wulbert, \emph{Extending {$n$}-convex functions}, Studia
  Math. \textbf{171} (2005), no.~2, 125--152. \MR{2183481 (2006h:26008)}

\bibitem{Pop34}
Tiberiu Popoviciu, \emph{Sur quelques propriétés des fonctions d'une ou de deux
  variables réelles}, Mathematica, Cluj \textbf{8} (1934), 1--85 (French).

\bibitem{Ral65}
Anthony Ralston, \emph{A first course in numerical analysis}, McGraw-Hill Book
  Co., New York, 1965. \MR{0191070 (32 \#8479)}

\bibitem{RobVar73}
A.~Wayne Roberts and Dale~E. Varberg, \emph{Convex functions}, Academic Press
  [A subsidiary of Harcourt Brace Jovanovich, publishers], New York-London,
  1973, Pure and Applied Mathematics, Vol. 57. \MR{0442824 (56 \#1201)}

\bibitem{Szeg75}
Gábor Szeg\H{o}, \emph{Orthogonal polynomials}, fourth ed., American
  Mathematical Society, Providence, R.I., 1975, American Mathematical Society,
  Colloquium Publications, Vol.~XXIII. \MR{0372517 (51 \#8724)}

\bibitem{WeiGL}
E.~W. Weisstein, \emph{Legendre-gauss quadrature}, from MathWorld-A Wolfram Web
  Resource.\\
  \texttt{http://mathworld.wolfram.com/Legendre-GaussQuadrature.html}.

\bibitem{WeiL}
\bysame, \emph{Lobatto quadrature}, from MathWorld-A Wolfram Web Resource.\\
  \texttt{http://mathworld.wolfram.com/LobattoQuadrature.html}.

\bibitem{WeiR}
\bysame, \emph{Radau quadrature}, from MathWorld-A Wolfram Web Resource.\\
  \texttt{http://mathworld.wolfram.com/RadauQuadrature.html}.

\bibitem{Was06JIPAM}
Szymon Wąsowicz, \emph{On error bounds for {G}auss-{L}egendre and {L}obatto
  quadrature rules}, JIPAM. J. Inequal. Pure Appl. Math. \textbf{7} (2006),
  no.~3, Article 84, 7 pp. (electronic). \MR{2257283 (2007e:41030)}

\bibitem{Was07JIPAM}
\bysame, \emph{Inequalities between the quadrature operators and error bounds
  of quadrature rules}, JIPAM. J. Inequal. Pure Appl. Math. \textbf{8} (2007),
  no.~2, Article 42, 8 pp. (electronic). \MR{2320615 (2008e:41020)}

\bibitem{Was07JMAA}
\bysame, \emph{Support-type properties of convex functions of higher order and
  {H}adamard-type inequalities}, J. Math. Anal. Appl. \textbf{332} (2007),
  no.~2, 1229--1241. \MR{2324332 (2008a:26028)}

\bibitem{Was08MIA}
\bysame, \emph{Hermite-{H}adamard-type inequalities in the approximate
  integration}, Math. Inequal. Appl. \textbf{11} (2008), no.~4, 693--700.

\bibitem{Was08JIPAM2}
\bysame, \emph{On quadrature rules, inequalities and error bounds}, JIPAM. J.
  Inequal. Pure Appl. Math. \textbf{9} (2008), no.~2, Article 36, 4 pp.
  (electronic).

\end{thebibliography}

\end{document}